\input harvmac
\baselineskip=.55truecm

\Title{\vbox{\hbox{alg-geom/9702002}}}
{\vbox{\centerline{Principal bundles on elliptic fibrations}}}
\vskip .05in
\centerline{\sl Ron Y. Donagi}
\vskip .2in
\centerline{\it Institute for Advanced Study}
\centerline{\it Princeton, NJ 08540, USA}
\centerline{\it and}
\centerline{\it Department of Mathematics,}
\centerline{\it University of Pennsylvania}
\centerline{\it Philadelphia, PA 19104, USA}

\vskip .2in
\centerline{ABSTRACT}
\vskip .2in
\noindent A central role in recent investigations of
the duality of F-theory and heterotic strings is played by 
the moduli of principal bundles, with various structure 
groups G, over an elliptically fibered Calabi-Yau manifold 
on which the heterotic theory is compactified. 
In this note we propose a simple 
algebro-geometric technique for studying the moduli spaces of 
principal G-bundles on an arbitrary variety X which is 
elliptically fibered over a base S: The moduli space 
itself is naturally fibered over a weighted projective 
base parametrizing spectral covers  $\tilde{S}$ of S, 
and the fibers are identified as translates of distinguished 
Pryms of these covers. In nice situations,  the generic Prym 
fiber is isogenous to the product of a finite group and an 
abelian subvariety of $Pic(\tilde{S})$.
\vskip .2in

\Date{\bf January 1997}

\def\C{{\cal C}}
\def\c{\underline{{\bf c}}}
\def\Hom{\rm{Hom}}
\def\Pic{\rm{Pic}}
\def\ad{{\bf ad}}
\def\Ad{\rm{Ad}}
\def\ker{\rm{ker}}
\def\mtx{{\cal M}^T_{X/S}}
\def\mte{{\cal M}^T_E}
\def\mtes{{\cal M}^T_{E_s}}
\def\mgx{{\cal M}^G_{X/S}}
\def\mge{{\cal M}^G_E}
\def\Aut{\rm{Aut}}
\def\psn{\overline{P_S/N}}
\def\pst{\overline{P_S/T}}
\def\ue{{\cal U}_E}
\def\ueb{\overline{\cal U}_E}
\def\ux{{\cal U}_{X/S}}
\def\uxb{\overline{\cal U}_{X/S}}
\def\O{{\cal O}}

\newsec{Introduction}

Moduli spaces of principal $G$-bundles on K3 surfaces, and on
 Calabi-Yaus of other dimensions, are basic ingredients of the 
compactification of heterotic strings, especially for semisimple 
structure groups $G$ contained in $E8 \times E8$ or in $Spin(32)/{\bf Z}2$. 
More recently, duality of the heterotic string with $F$-theory 
\ref\vafa{C. Vafa, `Evidence for F-theory', Nucl. Phys. 
B469(1996),403, HepTh 9602022; D. Morrison and C. Vafa, 
`Compactifications of F-theory on Calabi-Yau threefolds: 
I, II'', HepTh/9602114, HepTh/9603161.}  
has suggested the special importance of moduli spaces of 
$G$-bundles on elliptically fibered Calabi-Yaus. The purpose of this note
 is to propose a general technique for studying $G$-bundles on 
(not necessarily CY) elliptic fibrations $X \to S$, with 
arbitrary semisimple structure group $G$.

The idea is quite simple. To a principal $G$-bundle $P$ over $X$ we 
associate its {\it spectral data}. This consists of the {\it cameral 
cover} $\tilde{S} \to S$, a Galois cover  with  covering group the Weyl
 group $W$ of $G$; of a collection of multisection maps $v_{\lambda}:
 \tilde{S} \to X$, one for each character $\lambda \in \Lambda$ of
 the maximal torus $T \subset G$, subject to the condition of 
$W$-equivariance  (the {\it spectral covers} are the various images
 of $\tilde{S}$ in $X$); and of a $T$-bundle over $\tilde{S}$, subject 
to a certain twisted form of $W$-equivariance. When the cameral 
cover is reasonably nice, it determines the {\it distinguished Prym}
${\rm Prym}(\tilde{S}/S)$, an extension of an abelian variety by a 
finite group. The $G$-bundles on $X$ with a given cameral cover
 $\tilde{S}$ (and any consistent collection of maps to $X$) can then
 be parametrized by the distinguished Prym. An open subset of the moduli
 of $G$-bundles can thus be fibered over the parameter space for nice 
cameral covers, the fibers being the distinguished Pryms. When the 
cameral cover is not-so-nice, the description becomes less precise; but
 the whole construction has a local character, so bad behavior can be 
traced to specific singular loci. We avoid much of the difficulty by 
considering {\it regularized bundles}, or bundles with some additional 
structure, specified by a reduction of the structure group, over each 
$s \in S$, to an abelian subgroup which is the 
centralizer in $G$ of a regular element. (We call
 such a subgroup a regular centralizer.) Nice bundles have a unique 
regularization, special ones may have large families of regularizations,
 and others will have none, but can become regularized 
after some blowing up in the base $S$. 

The reason that regular centralizers are the right structure can be seen 
already from the behavior of $G$-bundles on a single elliptic curve $E$,
and already for $G=SL(2)$.
The structure group of ``most" semistable $G$-bundles on $E$ can be 
reduced to a maximal torus $T \subset G$. For $G=SL(2)$, there
is an essentially unique
exception: the non-trivial extension of $\O$ by itself (or the same 
tensored with one of the four spin bundles). The structure 
group in this case can still be reduced to an abelian subgroup, 
though not to $T$: it can be reduced to the group of upper triangular 
matrices with 1 or -1 on the diagonal, i.e. to
the centralizer of a non-zero nilpotent element in $SL(2)$.
As an abstract group, this is ${\bf Z}_2$ times the additive group 
$(\bf C,+)$. Now the family of bundles which admit a reduction to a 
given torus $T$ is parametrized by  ${\bf P}^1$, or more naturally by
$E$ modulo its involution. There are only 8 isomorphism classes of 
bundles with reduction to a nilpotent regular centralizer: the above 
four extensions of $L$ by $L$, where $L$ is a spin bundle, as well as 
the four trivial extensions $L \oplus L$. An important point is that all 
regular centralizers in $G$ fit nicely into one family, $\cal C$. 
For $SL(2)$, the base of 
this family is ${\bf P}^2 = ({\bf P}^1 \times {\bf P}^1)/ {\bf Z}_2$: 
the tori are parametrized by the complement of the diagonal ${\bf P}^1$, 
while the exotics sit over ${\bf P}^1$. 
Over this ${\bf P}^2$ lives the space 
parametrizing  $SL(2)$-bundles on $E$ together with a 
reduction of their structure 
group to {\it some} regular centralizer. This space looks like the threefold 
${\bf P}^1 \times {\bf P}^1 \times E$ divided by the obvious ${\bf Z}_2$, 
with the singularities blown up (yielding fibers of type $I^*_0$ over 
${\bf P}^1$) and the components of multiplicity $2$ discarded. Over this 
threefold there is a ``universal $\C$-bundle", from which our $G$-bundles
can be induced.
We will see that this picture generalizes to any $G$, and is the basic 
ingredient behind the reconstruction of a $G$-bundle from its spectral data.

Our construction is motivated by our previous work 
\ref\D{R. Donagi, Spectral covers, in:  Current topics in complex algebraic
 geometry, Math. Sci. Res. Inst. Publ. 28 (Berkeley, CA 1992/92), 65-86, 
AlgGeom 9505009 },
and in a sense is contained in it. The objects studied in \D \ may appear, 
at first sight, different than the ones that concern us here. There we were 
mostly interested in $K$-valued Higgs fields on $S$, where $K$ is a line 
(or vector) bundle on $S$. As we review in section 4, the parametrization
of these ``$K$-valued" objects is reduced to that of abstract, or unvalued,
 objects, and those are precisely the ones that come up in connection with 
$G$-bundles on $X \to S$. From this point of view, we may think of a 
$G$-bundle on $X$ as a Higgs bundle on $S$ whose ``values" are in the 
elliptic curves of the fibration.

The main construction, or rather its reduction to the construction in 
\D, is in section 4. Some basics on $G$-bundles on a single elliptic
 curve are gathered in section 2, while the behavior in families is 
discussed in section 3. There we also address the technical question 
of how to describe the family of all cameral covers. For this we need 
to understand the global properties of the moduli of $G$-bundles on a 
fixed elliptic curve, as well as the modular behavior seen when the
curve is varied. The former were worked out by Looijenga 
\ref\loij{E. Looijenga,  `Root systems and elliptic curves', Inv. 
Math. 38(1976),17-32 and `Invariant Theory for generalized root 
systems', Inv. Math. 61(1980),1-32}
and Bernstein and Shvartsman
\ref\ber{J. Bernstein and O.V.Shvartsman, `Chevalleys theorem for 
complex crystallographic Coxeter groups', Funct. Anal. Appl. 
12(1978),308-310},
and the latter, for all simply connected groups {\it except} 
$E8$, by Wirthmuller
\ref\wirt{K. Wirthmuller, `Root systems and Jacobi forms', 
Comp. Math. 82(1992),293-354.}.

It is a pleasure to thank Ed Witten, for asking the questions about 
moduli of $G$-bundles on elliptic fibrations which got me interested 
in the subject, and for drawing my attention to references \loij, \ber.
His joint  work with Friedman and Morgan 
\ref\fmw{R. Friedman, J. Morgan and E. Witten, `Vector bundles and 
F-theory', HepTh 9701162},
just posted to hep-th, has some overlap with this 
note. Roughly speaking, the emphasis in \fmw \ is on a description 
of the parameter space for spectral covers, while we focus
 on the fiber, which parametrizes bundles with a given cover.
 The theorem of Looijenga is recovered in full in \fmw, and
 there is also a discussion of Wirthmuller's work. 
Friedman, Morgan and Witten also 
obtain applications to the duality between F-theory
and the heterotic string. They present several beautiful descriptions 
of the moduli space of $G$-bundles on an elliptic curve. While we use 
deformations (of semisimple and semistable bundles, in the beginning of 
section 2) only to obtain a rather rough local picture of this moduli 
space, they base their main construction on the deformations of a
 ``minimally unstable" bundle. This gives them a global description,
 as well as the proof of Looijenga's theorem. In our approach, we 
take the results of Looijenga, Bernstein-Shvartsman, and Wirthmuller
 as given, and use them to describe the parameter space of the covers.
  The main novelty of our approach is the direct construction of the 
bundle corresponding to given spectral data, via the regular centralizers.

I am also grateful to Tony Pantev, for many discussions and for showing 
me an advance copy of another related work 
\ref\tony{M. Bershadsky, A. Johansen, T. Pantev and V. Sadov, `On
 four-dimensional compactifications of F-theory', hep-th 9701165}, 
to Eduard Looijenga who explained his works \loij \  to me and who told 
me about \wirt , and to Pierre Deligne for some very helpful comments 
on the manuscript. I have also benefitted greatly from conversations with
 J. Bernstein, R. Lazarsfeld, R. Livne, and V. Sadov on 
various aspects of the work described here. This work was partially 
supported by NSF grant DMS95-03249, a Lady Davis Fellowship from the 
Hebrew University, and grants from the  University  of Pennsylvania 
Research Foundation  and the Harmon Duncombe Foundation.     

\newsec{G-bundles on elliptic curves}

Let $G$ be a connected, simply-connected
 complex semisimple group, $T$ its maximal torus,
 $E$ an elliptic curve, i.e. a non singular curve of genus 1 with a marked 
point $0\in E$.  The moduli space $M^T_E$ of 
degree-0 semistable $T$-bundles on $E$ 
is $\Hom (\Lambda, \check{E})$, where $\Lambda := \Hom (T, {\bf C})$ is the 
lattice of characters of $G$, and $\check{E}:=\Pic^0 E$ is the dual elliptic
 curve, which is naturally identified with $E$.
 This moduli space is an abelian
 variety, in fact it is (non-naturally) isomorphic to $E^r$,
 where $r$ is the
 rank of $G$.  It comes with a natural action of the Weyl group $W$, 
as well as a natural polarization (cf.\loij)
 whose degree is
 the discriminant of $\Lambda$.  Over $E \times M^T_E$ there
 is a universal
 $T$-bundle. In case $T={\bf C}^*$, $M^T_E$ is just the dual
 elliptic curve
 $\check{E}$, and the universal bundle is the Poincare bundle.
  Via the
 unnatural identification of the general $M^T_E$ with $E^r$,
 the universal
 bundle becomes the sum of $r$ Poincare bundles pulled back
 from the $r$
 factors.  In this abelian situation there is actually
 a moduli space of all
 (not necessarily semistable) $T$-bundles: just replace
 $\check{E}=\Pic^0
 E$ throughout by $\Pic \  E$.

The moduli space $M^G_E$ of semistable $G$-bundles on $E$ 
is the quotient
 ${M^T_E}/{W}$, where the Weyl group $W$ acts through
 $\Lambda$.  This moduli
 space parametrizes not isomorphism classes but $s$-equivalence
 classes of
 semistable $G$-bundles.  There is an open set where the two
 notions coincide,
 but over smaller strata there will be a finite number of 
isomorphism classes
 represented by each point in moduli. 
 The set of isomorphism classes maps
 naturally to $M^G_E$, and this map behaves
 nicely in families.  (This
 classifying map is, of course,  part of the definition of
a coarse moduli space.) 
This map can be given as follows: 
any semistable principal $G$-bundle 
$P$ on $E$ has a semistable
 reduction of its structure group to the Borel
 subgroup $B\subset G$, in other words, it is associated to some 
semisimple $B$-bundle $P_B$. 
 The quotient map $B\rightarrow T$ then gives an
 associated $T$-bundle 
$P_T$ which, up to the action of $W$, is 
independent of the choice of Borel 
reduction.  The existence of a Borel reduction
 can be seen, e.g., by finding
 a flat holomorphic connection on $P$
 (\ref\ginz{V. Baranovsky and V. 
Ginzburg, Conjugacy classes in loop groups 
and G-bundles on elliptic curves
 AlgGeom 9607008}) and noting that the 
holonomy is abelian, as image of
 $\pi_1(E)$. (The entire discussion extends easily to 
non simply connected groups, as long as we restrict attention 
to those bundles which are liftable to
the simply connected covering group.)

We can develop a feel for this moduli 
space by considering the deformation 
theory of semisimple bundles, i.e. ones 
associated to $T$-bundles.  Let $P_T$ 
be a principal $T$-bundle on $E$,
 and $P$ the associated (semisimple)
 $G$-bundle.  Deformations of $P$ are 
unobstructed, parametrized by $H^0(\ad P)^*$.  Now 
$$\ad (P)=\ad (P_T)\oplus(\oplus_\alpha L_\alpha),$$ 
where the sum is over the roots $\alpha \in R$ of $G$,
 and $L_\alpha$ is 
the line bundle associated to $P_T$ by the root $\alpha$.
  For generic
 $P$ each $L_\alpha$ is a non-trivial line bundle of
 degree 0 on $E$, 
so $H^0 (\ad P)=H^0 (\ad P_T)$.  This means that all
 deformations of such
 generic $P$ come from deformations of $P_T$.
But for non-generic semisimple $P$, there can be other 
deformations. The possible types of non-generic 
semisimple $P$ are indexed by subroot systems:
 An arbitrary $P_T$ is given
 by a homomorphism $p:\Lambda \rightarrow E$,
 which determines a sub root
 system $R':=R\cap \ker (p)$, generating 
a sublattice ${\Lambda}'
 \subset \ker (p) \subset \Lambda$.  The structure
 group of this $P_T$ then
 reduces to $\bar{T}:=\Hom ({\Lambda}/{{\Lambda}'},
 \, {\bf C}^*)\subset T$.
  Let $G'$ be the centralizer $Z_G(\bar{T})$, with Lie
 algebra $\bf{g}'$.  The $\alpha$ for which $L_\alpha$ 
is trivial are precisely the roots of $\bf{g}'$, so the
 deformations are now parametrized by the
 dual $\bf{g}'^*$, and the general
 deformation has structure group $G'$.  Since
 $G'$ is the automorphism group
 of $P$, it acts on the deformation space, and the
 isomorphism type is fixed
 along $G'$-orbits.  The semisimple locus looks near
 $P_T$ like  ${\bf{t}}$, and modulo $G'$ this becomes
 ${\bf{t}}/{{W}'}$ (${\bf t}$ is the Cartan
 subalgebra, ${W}'$ is
the Weyl group of ${R}'$), while the new transversal
 deformations
 look like the nilpotent cone in $\bf{g}'$.
Modulo $G'$, this gives just a finite set of types, 
with a non-separated topology.  
In case $G=SL(n)$ we recover
 Atiyah's description 
\ref\At{M. Atiyah, Vector bundles over an elliptic
 curve, Proc.
 LMS VII(1957),414-452.}
of vector bundles on elliptic curves as sums of
 indecomposables. 
 The possible root systems ${R}'$ correspond in this case to
 faces of the Weyl
 chamber, or equivalently to parabolic subgroups containing $B$.  
For other groups there
 may be additional possibilites; e.g. for $G_2$ there
 are, in addition
 to the parabolic types $0,\, A^{\rm{short}}_1,\,
 A^{\rm{long}}_1$, and
 $G_2$, also the two possibilities $A^{\rm{long}}_2$ and
 $A^{\rm{long}}_1
 \times A^{\rm{short}}_1$, which arise when the image 
of  $p:\Lambda \rightarrow E$ happens to consist of 
points of order two. But in any case, the
 nilpotent cone breaks
 into a finite number of orbits.  The closed orbit
 is 0, corresponding
 to the semisimple bundles.  At the opposite end,
 there is a unique
 dense orbit, corresponding to regular bundles.
  For such bundles
 $P, \, H^0(\ad P)$ is again of the smallest
 possible dimension,
 namely $r$.  It is the centralizer in $\bf{g}'$
 of a regular nilpotent
 element.

There is a natural, ``universal", space $\ue$
 parametrizing $G$-bundles
 on $E$ together with a trivialization (frame)
 over $0 \in E$ and a
 reduction to a regular centralizer.  To describe it, 
we start with the quotient ${G}/{T}$ which 
 parametrizes pairs $T' \subset B'$ consisting of a
 maximal torus $T'$ and a Borel subgroup containing it,
 or equivalently a torus with a choice of chamber.  The
 Weyl group $W$ acts on ${G}/{T}$.  The quotient ${G}/{N}$
 (where $N:=N_G(T)$ is the normalizer of $T$ in $G$)
 parametrizes tori $T'$.  Equivalently, ${G}/{T}$ and ${G}/{N}$
 can be described in terms of Cartan subalgebras
 $\bf{t'} \subset \bf g$  and Borels ${\bf b' \subset \bf g}$.
  In [D] we introduced the parameter spaces $\overline{{G}/{N}}$
 and $\overline{{G}/{T}}$, parametrizing all regular centralizers
 $C \subset G$, respectively pairs $C \subset
 B$. By definition, $\overline{{G}/{N}}$ is an open subset of 
the closure of ${G}/{N}$ in ${\rm Grass}(r, \bf g)$. 
 (The closure itself contains some abelian subalgebras
 which are not regular centralizers and therefore need to be discarded.)

Over $\overline{G/N}$ there is the vector bundle $\c \to
 \overline{G/N}$ of rank $r$ abelian subalgebras of $\bf g$,
 and the corresponding abelian group scheme $\C \to
 \overline{{G}/{N}}$ whose fibers are the regular 
centralizers in $G$. Let $\varpi$ be the projection 
$\overline{{G}/{N}} \times E \to \overline{{G}/{N}}$. 
Then $\varpi^*\C$ is an abelian group scheme
 over $\overline{{G}/{N}} \times E$. The cohomology
 sheaf  $R^1\varpi_*\varpi^*\C$ can be represented by an
 analytic or algebraic space $u': \ue ' \to \overline{G/N}$
 whose fiber over $C \in \overline{G/N}$ is
 $H^1(E, C({\cal O}_E))$, the moduli space of
 $C$-bundles on $E$. We will usually restrict attention to
 the subfamily  $u: \ue \to \overline{G/N}$ parametrizing
 semistable $C$-bundles. Over $\ue$ there is a universal principal
 $\C$-bundle (more precisely, a $u^*\C$-bundle) $P^C_{\ue} \to \ue
 \times E$. Since $u^*\C$ is a subgroup scheme of the trivial group
 scheme $G$ over $\ue \times E$, we get an associated principal
 $G$-bundle $P_{\ue} \to \ue \times E$.
The bundles  $P^C_{\ue}$ and  $P_{\ue}$ are uniquely characterized
 by the properties:

\noindent $*$ The restriction of $P^C_{\ue}$ to $\ue \times \{ 0 \}$
 is the trivial $\C$-bundle.

\noindent $*$ the restriction of  $P^C_{\ue}$ to $\{ x \} \times E$
 (where $x:=(C,p) \in \ue$, $C$ a regular centralizer, $p$ the isomorphism
 class of a  $C$-bundle over $E$) is a $C$-bundle on $E$ in the class
  $p$, and the restriction of $P_{\ue}$ is the associated $G$-bundle.

We could rephrase this as saying that $\ue$ is a {\it fine} moduli space
 for the data it parametrizes: a $G$-bundle together with a
 trivialization over $0 \in E$ and a semistable reduction to 
a regular centralizer. The basic reason for existence of the universal
 family is that the objects parametrized have no automorphisms: we
 killed them by fixing the trivialization.

Over any one stratum of $\overline{G/N}$ it is easy to describe $\ue$
 and $P^C_{\ue}$. For instance, over the open stratum $G/N$, we start
 with the Poincare $T$-bundle over $\mte \times E$, cross this with
 $G/T$, and divide by $W$ which acts on both $G/T$ and $\mte$.  

A related object which we shall need is the quotient $\ueb :=
 (\overline{G/T} \times \mte )/W$.  There is a natural morphism
 $f: \ue \to \ueb$ which, over each $C \in \overline{G/N}$, maps
  $H^1(E, C({\cal O}_E))$ to its compact part
  $H^1(E, C'({\cal O}_E))$, where $C'$ is the quotient of
 $C$ by its unipotent part, $C' := C / (C \cap [B,B])$ for a
 Borel $B \supset C$. Between fibers over points of
 the open stratum $G/N$, this map is surjective (in fact, an isomorphism);
 but over the whole base its image is the constructible set 
$${\ueb} ':=\{ (C,t) \in  \overline{G/T} \times \mte \ \ \vert \ \ 
{\rm Stab}_W(C) \subset {\rm Stab}_W(t)\}/W .$$
As described in the introduction, for $G=SL(2)$ we have 
$\overline{G/T} = {\bf P}^1 \times {\bf P}^1$ and  
$\overline{G/N}={\bf P}^2$. The fiber of either $\ue$ or ${\ueb}'$ 
over a point not in the diagonal is $E$. Over points of the diagonal, 
the fibers in $\ue,\ {\ueb}'$ and $\ueb$ respectively are: 
four lines, four points, and ${\bf P}^1$. 
We may think of \ ${\ueb}'$ as parametrizing isomorphism classes of
 $G$-bundles having a semistable reduction to a regular centralizer
 $C$, together with a trivialization over $0 \in E$. The additional
 data in the fiber of $\ue$ over $\ueb$ chooses such a $C$-structure.
 Thus $f$ is the forgetful map sending a $C$-struture to its associated
 $G$-structure. The point is that the dimension of the normalizer
 $N_G(C)$ can be greater than $r$, so there can be a non-trivial
 family of $C$-bundles whose associated $G$-bundles are all isomorphic.

\newsec{Families of moduli spaces}

Let $\pi : X \to S$ be an elliptic fibration with a section, with
 non-singular $X,S$. We want to put the basic objects from the previous
 section into families parametrized by $S$. It will be
 convenient to work instead with a (singular)  Weierstrass model
  $\overline{X} \to S$, given by an equation of the form
 $y^2 = x^3 +b_2x+b_3$, where $b_i$ is a section of
 $L^{\otimes 2i}$ and $L$ is $\pi_*$ of the relative canonical
bundle $K_{\pi}$. 

Since $\pi$ has a section, we can identify $\Pic ^0$ of a smooth fiber
 $E_s := \pi^{-1}(s)$ with $E_s$. Globally, the relative Jacobian
 $\Pic ^0 (\bar{X}/S)$ is the complement in $\bar{X}$ 
of the locus of singular
 points of fibers. (By definition, this singular locus includes any
 fiber component of multiplicity $>1$.) The reason for this is that
 a section of  $\pi : \bar{X} \to S$ must have intersection number $1$
 with each fiber, so it cannot pass through a singular point. This
 leads immediately to identification of $\mtx$ with the $r$-th
 cartesian power of $\bar{X} - {\rm Sing}(\bar{X})$ over the base $S$.

A satisfactory construction of $\mgx$ is somewhat more delicate.
 Looijenga shows \loij \ (see also  \ber ) that 
for simply connected group
 $G$ and for fixed, non-singular $E$, $\mge$ is a weighted
 projective space. The weights are the Dynkin indices of the
 dual Dynkin diagram of $G$ (i.e. the coefficients of the 
highest short root when expressed as a sum of simple roots, 
plus the coefficient 1 for the affine root). 
In order to be able to describe
 the relative object $\mgx$, we need a way of relating nearby
 fibers. A general construction of a flat connection on such
 families was carried out by Saito
\ref\sai{K. Saito, Extended affine root systems, Publ. RIMS
 Kyoto,  I: 21(1985), 75-179, and II: 26(1990),15-78.}
but more immediately useful results are in \wirt . Wirthmuller
 takes the base $S$ to be a modular curve, so that Jacobi forms
 give sections of powers of the ``theta" bundle (giving the
 polarization) on $\mtx$, and we are looking for the
 $W$-invariants among them. These are now bi-graded, by what
 he calls ``weight" (as modular form) and index. The ``indices"
 correspond to the weights in Looijenga's weighted projective
 space, and also to the power of the theta bundle in which the
 section lives.  I will refer to Wirthmuller's ``weights" as
 degrees, to avoid confusion with Looijenga's weights. These
 degrees correspond to powers of $L$, the Weierstrass line
 bundle. One of them always turns out to be $0$, the others are
 the degrees of the basic $G$-invariant polynomials. It is not
 clear a-priori why there should be a natural way to pair these
 two sequences, of indices and degrees. The end result is that
 $\mge$, for simple $G$ other than $E8$, can be identified with
 the quotient $(\oplus_{i=0}^r L^{d_i})/{\bf C}^*$, where
 $d_0=0$, $\{ d_i \}_{i=1}^r$ are the degrees of the invariant
 polynomials, and ${\bf C}^*$ acts with weights equal to the
 Dynkin indices for the dual Kac-Moody algebra.

At least for some groups, these results are elementary. For
 $SL(r+1)$, all indices equal $1$, while the degrees are
 $0,2,3,4,...r+1$. So $\mge$ is the ordinary projective space,
 projectivization of $\O \oplus L^2 \oplus ... \oplus L^{r+1}$.
 The fact that $\mge$ is a projective space follows directly from
 Abel-Jacobi: $\mte$ is the locus of $(r+1)$-tuples of points in
 $E$ which add up to $0$, $W$ is the symmetric group, so the
 quotient is the variety of effective divisors in the linear
 system $r+1$ times the origin, a projective space of dimension
 $r$. But we can really identify this space: for $r=1$ it is the
 ${\bf P}^1$ with coordinate x, of which the Weierstrass equation
 exhibits $E$ as a double cover. A basis of sections of $\O(1)$ is
 given by $1$ (of degree 0) and $x$ (of degree $2$). For $r=2$ we
 get the ${\bf P}^2$ with functions $1,x,y$, etc. Similar or easier
 arguments show that for $G$ either $SO(2r+1)$ or $Sp(r)$, the moduli
 space $\mge$ is the $r$-th symmetric product of the basic ${\bf P}^1$,
 so it can be identified with the projectivization (all weights equal 1)
 of $\O + L^2 + L^4 + ... +L^{2r}$. For $Sp(r)$, this is as it should be:
 the Dynkin diagram of the Kac-Moody
 algebra $\check{(C_r)}=D_{r+1}^{(2)}$
 has all indices equal to 1. But the Dynkin diagram of the Kac-Moody
 algebra $\check{(B_r)}=A_{2r-1}^{(2)}$ has its three extreme indices
 equal to $1$, while the others equal $2$, so we expect a weighted
 projective space of weights $(1,1,1,2,...,2)$. And that is exactly
 what we get, if we replace $SO(2r+1)$ by the simply connected
 $Spin(2r+1)$: the double covering group yields a $4$-sheeted covering
 moduli space (the elliptic curve has $4$ two-torsion points), which
 can be identified as the unique ${\bf Z}/2 \times {\bf Z}/2$ cover
 each of whose three ${\bf Z}/2$ quotients is branched along two of
 the three hyperplanes in ${\bf P}^r = {\rm Sym}^r {\bf P}^1$ 
corresponding to $r$-tuples containing one of the three roots of 
the Weierstrass equation.

I have not worked out the missing case of $E8$, nor the remaining non
 simply connected groups. Presumably, even if $\mge$ does not always
 turn out to be a weighted projective space, it can still be described
 by some universal construction in terms of the Weierstrass line bundle
 $L$. In particular, there is then a natural way to extend the family
 of quotients $\mtes /W$ to a locally trivial family over the entire
 base $S$ (including the discriminant) of any Weierstrass family. The
 most familiar instance of this is when $G=SL(2)$: there the resulting
 $\mgx$ is the ${\bf P}^1$-bundle over $S$ obtained as projectivization
 of $\O \oplus L$, which is the quotient of the Weierstrass model by
 its natural involution fixing the $0$-section. 

We also need to note that the spaces $\ue$ and $\ueb$ parametrizing
 $\C$- and $G$- bundles likewise extend in families to form objects
 $\ux$ and $\uxb$ respectively. For $\uxb$ we take $( \overline{G/T}
 \times \mtx )/W$. Near $s$ corresponding to non-singular fiber $E_s$,
 this looks like a bundle over $S$ with fibers $\ueb$. At singular
 $s$, we have already removed the singularities of $\mtx$, so we get
 a non-singular, but somewhat smaller, object. Likewise, we restrict
 $\ux$ to the open subset which maps to $\uxb$. Its points then
 parametrize $\C$-bundles on the $E_s$, and there is again a universal
 $\C$-bundle on $\ux \times _S X$. (It is possible that for those
 groups for which  Wirthmuller-type constructions work we could extend
 these universal objects to some of the singular loci, but we have not
 pursued this, and the construction presented below seems to avoid the
 issue.)

\newsec{Spectral parametrization of bundles}

In this section we describe the equivalence between regularized
 $G$-bundles on an elliptic fibration with a section and the
 corresponding spectral data. We fix an elliptic fibration
 $\pi: X \to S$ with a ``zero" section $\sigma: S \to X$, where
 $X,S$ are smooth. Let $E_s := \pi^{-1}(s)$ denote a fiber.
 Given a principal $G$-bundle $P$ on $X$, let
 $P_S := P_{\vert \sigma (S)}$ be the restriction
 to the zero section.  There are associated bundles
 of groups $\Ad P, \Ad P_S$ over $X,S$ respectively.
 The sheaf of automorphisms along the fibers
 $E_s$, $\Aut _S (P) := \pi_* Ad P$, can be identified
 as a subsheaf of $\Ad P_S$.   There is also an associated
 bundle $\overline{P_S/N}$, whose fiber over $s \in S$ is the
 family of regular centralizers associated to the fiber $P_s$
 over $\sigma (s)$.

A section  $c: S \to \psn$
determines an abelian group scheme $\C \to S$, a subgroup scheme
 of the group scheme $\Ad P_S$. A $\pi^* \C$-torser $P^{\C}$ on $X$ 
then induces a $G$-bundle $P=P^{\C} \times _{\C}  \Ad P_S .$
By a {\it regularized $G$-bundle}
 we mean such a triple $\{ P, \ c:S \to \overline{P_S/N}, \ P^{\C} \},$ 
or equivalently a reduction of the structure of $P$ to 
a group scheme $\C$ of regular centralizers.
This group subscheme $\C \subset \Ad P_S$ is contained in 
 $\Aut _S (P)$, and extends naturally to a group 
subscheme of $\Ad P$.
A key point is that an everywhere regular, semisimple and semistable
 bundle $P$ (i.e. one whose restriction $P_s$ to each $E_s$ has these 
 properties) has a unique regularization, with $\C =\Aut _S (P)$.
 But there are other regularized bundles, whose underlying
 $G$-bundles may not be everywhere, or even anywhere, regular. Yet
 these bundles too can be parametrized by their spectral data.
 (Being regularized means that we have chosen a reduction to a regular 
 centralizer subgroup of the automorphisms along each $E_s$, not necessarily 
 that those automorphisms form a regular centralizer themselves.)
If $P$ is known to be regular, semisimple, and semistable only for generic 
$s \in S$ then a regularization is still unique if it exists; in general 
though, we may have to blow up the base $S$ to find  a regularization.

A cameral cover of $S$ is a $W$-Galois cover $\tilde{S} \to S$ which is
 modelled on $\overline{G/T} \to \overline{G/N}$. (``modelled on" means
 ``obtained locally as pullback via maps of $S$ to the base".) Recall
 that $\overline{G/N}$ parametrizes regular centralizers in $G$, while
 $\overline{G/T}$ parametrizes pairs of a regular centralizer and a
 Borel containing it.

A regularized bundle $\{ P, c , P^{\C} \}$ determines the following data: 

\noindent (1) A cameral cover $\tilde{S} \to S$.

\noindent (2) A $W$-equivariant morphism  $v: \tilde{S} \to \mtx$.
 (Equivalently, a $W$-invariant morphism $v': \Lambda \times \tilde{S} \to
 \Pic ^0(X/S)$. By way of terminology, we refer to the
 image of $v$ as the {\it universal spectral cover}. The
 various other spectral covers are the images under $v'$
 of $\lambda \times \tilde{S}$, for $\lambda \in \Lambda$.)

\noindent (3) A homomorphism $\ell : \Lambda \to \Pic
 \ \tilde{S}$ (or equivalently, a $T$-bundle on
 $\tilde{S}$) satisfying the twisted $W$-equivariance
 property of \D .

These are obtained as follows. 

\noindent (1) The cameral cover $\tilde{S}$ is the cover of
 $S$ induced via $c$ from the cover $\overline{P_S/T} \to
 \overline{P_S/N}$, which  indeed looks locally like
  $\overline{G/T} \to \overline{G/N}$.

\noindent (2) A point $\tilde{s} \in \tilde{S}$ above
 $s \in S$ corresponds to a choice of Borel in $P_s$ containing
 the regular centralizer $\C _s \subset \Aut _s (P_s)$. This
 extends uniquely to a subbundle of Borels in $P_{\vert E_s}$.
 Via the natural quotient map $B \to B/[B,B] = T$, there is an
 associated $T$-bundle, identified with a point of $\mtes$. 
 When $\tilde{S}$ is reduced,  we get the map $v$ by letting
 $s$ and $\tilde{s}$ vary, and the $W$-equivariance
 holds, since it does so fiber-by-fiber. In the general case,
 the $\C$-bundle $P^{\C}$ on $X$ together with the inclusion of
 $\C$ into the universal Borel bundle $\cal B$ over $\tilde{S}$ 
 induce on $\tilde{S} \times _S X$ a $\cal B$-bundle, hence a
 quotient $T$-bundle.  Our morphism $v$ is the classifying map
 for this bundle.

\noindent (3) Above $\tilde{S}$ we have a bundle of Borels; the
 $T$-bundle is associated to it as above. We will discuss the shifted
 $W$-equivariance below.

The main point of this note is that a regularized $G$-bundle can be
 reconstructed from its spectral data. To see this, it will be
 convenient to introduce an intermediate object, the {\it principal
 G-Higgs bundle}, cf. \D. This is simply a pair $(P_S, \C)$ consisting of a
 principal $G$-bundle $P_S$ on $S$ together with a family $\C$ of
 regular centralizers in $\Ad \ P_S$. (To avoid worrying about how
 the different types of centralizers fit together, we can instead
  consider the {\it vector} subbundle $\c$ of regular centralizers
 in the bundle $\ad \ P_S$ of Lie algebras. This , of course, is
 equivalent data.)  Now as above, to a principal $G$-Higgs bundle 
$(P_S, \C)$ corresponds
 an abstract cameral cover $\tilde{S} \to S$ and a homomorphism
 $\ell : \Lambda \to \Pic \ \tilde{S}$. (The constructions of items
 (1) and (3) above used only the available data on $S$.) As $W$ acts
 on both sides, we can consider the subgroup of $W$-equivariant
 homomorphisms. In \D   \  this was called the {\it distinguished Prym}
 of $\tilde{S}/S$. In case $S, \tilde{S}$ are non-singular, this is a
 finite group times an abelian subvariety
 of $\Hom (\Lambda, \Pic (\tilde{S}))$.

The main observation in \D \ is that the family of principal Higgs
 bundles with a given cameral cover $\tilde{S}$ is, if non-empty,
 parametrized by a translate of the distinguished Prym. The exact
 point by which we need to translate will not be crucial for us 
here. It is described in sections 5.2 and 5.3 of \D \ as the sum 
of a cohomological shift term depending only on the group $G$, and
 further twist terms coming from the fixed divisors for the action
 of $W$ on $\tilde{S}$. An analogous and more familiar situation
 applies in case $G=SL(n)$ when we replace the cameral cover by
 the degree-$n$ spectral cover $\pi : \bar{S} \to S$ ($\bar{S}$
 is the image   under $v'$ of $\lambda_1 \times \tilde{S}$, where
 $\lambda_1$ is the first fundamental weight.) The twisting along
 the ramification then corresponds to the relative Todd class
 (=half the ramification) which enters into the
 Grothendieck-Riemann-Roch formula: In order for $\pi_*(L)$ to have
 determinant $0$, $c_1(L)$ needs to be the
 class of half the ramification.

(The original purpose of \D \ was to describe the fibers of the
 Hitchin map, from moduli of $K$-valued Higgs bundles on a variety
 $S$, to $K$-valued spectral data. Here $K$ can be the canonical
 bundle of $S$, as in Hitchin's original work \ref\hit{N. Hitchin,
  Stable bundles and integrable systems, Duke Math. J.
 54(1987),91-114.}, but it can also be an arbitrary line bundle,
 or (with some additional symmetry conditions imposed) even a
 vector bundle on $S$ as in Simpson's works. It turned out to be
 convenient to separate the problem into considerations of
 ``eigenvectors" and ``eigenvalues": we introduce the somewhat
 abstract ``principal G-Higgs bundles" and show that they
 correspond to abstract spectral data (this is the eigenvector
 aspect); the $K$-valued versions are then recovered by adding
 a $K$-valued ``Higgs field" $\phi$, a section of $\ad (P_S)
 \otimes K$, on the one side, and an ``eigenvalue map" 
 $v: \tilde{S} \to K$ on the other. This latter map is of
 course analogous to our datum (2). We are thus led to think
 of a $G$-bundle on an elliptic fibration as a sort of Higgs
 bundle on the base taking values in the fibration instead of
 in a line (or vector) bundle.)

Returning to our situation, from the spectral data (1)-(3) we
 thus retrieve the principal $G$-Higgs bundle  $(P_S, \C)$
 together with a morphism $v$ from $\tilde{S}$ (which is
 determined by $(P_S, \C)$) to $\mtx$, commuting with the
 projections to $S$. It remains to recover the original
 regularized $G$-bundle $\{ P, c , P^{\C} \}$ (on $X$) 
 from this data. This goes as follows: the situation is 
 essentially rigid, so we can reduce to the case that $S$ is a 
 point, i.e. to bundles on one elliptic curve. For a given
 regular centralizer $C$, this then reduces to the straightforward
 verification that the canonical map:
 $$ {\cal M}^C_E \to \rm{Mor}_W((G/B)^C, \mte ) $$
 is an isomorphism. Here $(G/B)^C$ is the subscheme of $G/B$
 parametrizing Borels through $C$ (it is finite of length equal
 to the cardinality of $W$, and is not reduced except when $C$
 is a torus), and $\rm{Mor}_W$ is the (group of) $W$-equivariant
 morphisms.

 Working globally over $S$, the principal $G$-Higgs bundle  $(P_S, \C)$
 determines an (``eigenvector") map $c: S \to \psn$ and its lift
 $\tilde{c}: \tilde{S} \to \pst$, which together with the ``eigenvalue"
 map $v: \tilde{S} \to \mtx$ sends $\tilde{S}$ to the fiber
 product $\pst \times _S \mtx$. This map is $W$-equivariant,
 so it descends to give a section $\bar{\alpha}: S \to \uxb ^{P_S}:=
 (\pst \times _S \mtx )/W$. This last object is a version
 of our previous $\uxb$ which is twisted by the bundle $P_S$:
 it is isomorphic to $\uxb$  over any open subset of $S$ over
 which $P_S$ is trivial. Likewise, we have the twisted version 
$\ux^{P_S}$ of $\ux$. Now, the data needed to lift our section 
$\bar{\alpha}$ to a section $\alpha$ of $\ux^{P_S}$ is precisely given by 
$v: \tilde{S} \to \mtx$. (In particular, existence of $v$ implies that the image of  $\bar{\alpha}$ is contained in ${\uxb}'$.) Now on
 $\ux \times _S X$ we have the universal $\C$-bundle trivialized
 along $S$. Twisting by $P_S$, we get  a universal $\C$-bundle
 on $\ux ^{P_S} \times _S X$ , but instead of a trivialization
 we now get an identification of its restriction to $S$ with
 the universal $\C$-subbundle of $P_S$. Pulling back via the
 section $\alpha$ to $X=S \times _S X$ gives a $\C$-bundle
 $P^{\C}$, whose associated $G$-bundle is the original $P$.

$\underline{\rm Remark 1}$ We have emphasized the cameral covers,
 since they are, in our view, the most basic objects in the picture.
 But in order to parametrize entire components of the moduli space of
 $G$-bundles on $X$, it is necessary to allow the cover to vary. For
 this purpose, it is more convenient to consider instead the universal
 spectral covers. These are all obtained by pulling back one object,
 the $W$-cover $\mtx \to \mgx$, by arbitrary sections $S \to \mgx$. 
As we saw in the previous section, $\mgx$ is, 
for most groups, a bundle of
 weighted projective spaces over S. So in a connected component
 (obtained by fixing the numerical invariants of the map), these
 maps are specified by the (weighted projectivization 
of) the space of sections
 of an appropriate vector bundle over $S$. In ``nice" situations, a
 generic map of this kind will determine a non-singular universal
 spectral cover, and therefore also a unique cameral cover.

$\underline{\rm Remark 2}$ An alternative approach to reconstruction
 of a $G$-bundle on $X$ from spectral data might be based on
 application of the equivalence of \D \ directly to the 
principal $G$-Higgs bundle $(P, \pi^*\C )$ on $X$. One then 
obtains a cameral cover $\tilde{X} \to X$ and a point 
$\ell _X$ in a translate of its distinguished Prym. 
It is interesting to compare this to the spectral data 
$(\tilde{S}, v, \ell _S)$ which we used above. The cover  
$\tilde{X} \to X$ is clearly the pullback via 
$\pi$ of  $\tilde{S} \to S$, but $\pi ^* \ell _S$ is
{\it not} the same as $\ell _X$: the former is trivial along 
fibers of $\pi$, while the behavior of the latter along fibers 
is equivalent to the additional datum $v: \tilde{S} \to \mtx$. 
Our original approach has a clear advantage over this 
alternative in situations where $X$ is interesting (say a 
$K3$ or $CY_3$) while $S$ is a much simpler object 
(${\bf P}^1, {\bf F}_n$).

\listrefs
\end